\documentclass[11pt]{article}
\usepackage{bbm}
\usepackage{mathrsfs}
\usepackage{amsmath,amsthm,enumerate}
\usepackage{amsfonts,amssymb,color}
\usepackage{dsfont}
\usepackage{curves}
\usepackage{mathrsfs}
\usepackage{pifont}
\usepackage{graphicx}
\usepackage{float}
\usepackage{hyperref}
\usepackage{cite}

\newtheorem{theorem}{Theorem}[section]
\newtheorem{lemma}{Lemma}[section]
\newtheorem{corollary}{Corollary}[section]

\newtheorem{remark}{Remark}[section]

\numberwithin{equation}{section}
\numberwithin{figure}{section}
\usepackage{amsmath}

\date{}

\def\emptyset{\mbox{{\rm \O}}}
\def\bar{\overline}
\DeclareMathOperator{\diag}{diag}

\begin{document}
	\title{\textbf{ Which $L$-cospectral graphs\\ have same degree sequences}\footnote{ Emails: yejiachang12@163.com (J. Ye).}}
	
	\author{\small   Jiachang Ye$^{1}$   \\ \small $^1$ School of Mathematical Sciences, Xiamen University,  Xiamen, 361005, P.R. China}

	\maketitle
	
	\begin{abstract}
Let $\lambda_{i}(G)$ be the $i$-th largest Laplacian  eigenvalues of graph $G$, where $1\le i\le |V(G)|$. Liu,  Yuan,  You and Chen [{\em Discrete Math.}, 341 (2018) 2969--2976] raised the problem for ``Which cospectral graphs have same degree sequences".  In this paper, let $W_3$ and $W_5$ be the two graphs as shown in Fig. \ref{fig2} and let $G$ be a connected  graph with $n\ge 18$ vertices. We shall show that:
\begin{itemize}
	\item[(1)]  If $\lambda_{2}(G)<5<n-1<\lambda_{1}(G)$, $\lambda_{1}(G) \notin \{\lambda_{1}(W_3),\lambda_{1}(W_5)\}$ and $H$ is Laplacian cospectral with $G$, then  $H$ must  have the  same degree sequence with $G$; \item[(2)] If $\lambda_2(G)\le 4.7<n-2< \lambda_1(G)$, and $H$ is   Laplacian  cospectral with $G$,  then  $H$ must  have the  same degree sequence with $G$.   \end{itemize}
  The former result  easily leads to the unique theorem result of [{\em Discrete Math.}, 308 (2008) 4267--4271], that is: Every multi-fan graph $K_1\vee (P_{l_1}\cup P_{l_1}\cup\cdots \cup P_{l_t})$ is determined by the Laplacian spectrum. Moreover, it can also deduce a new conclusion: $K_1\vee (P_{l_1}\cup P_{l_1}\cup\cdots \cup P_{l_t}\cup C_{s_1}\cup C_{s_2}\cup\cdots \cup C_{s_k})$ $(t\ge 1, k\ge 1)$ is determined by the Laplacian spectrum if  the graph order $n\ge 18$ and each $s_i$ $(i=1,2,\ldots, k)$ is odd.
		
		\begin{flushleft}
		\end{flushleft}
	\end{abstract}

\section{Introduction}
	In this paper, $G =(V,E)$ is a simple undirected graph with $n$ vertices. For $ v\in V(G) $,  denote $d_G(v)$ and $N_G(v)$ as the degree and neighbourhood set of  $v$ in the graph $G$. In the sequel, the degrees of $G$ are enumerated in non-increasing order, i.e., $d_{1}(G)\ge d_{2}(G)\ge ...\ge d_{n}(G)$, where $d_{G}(v_{i})=d_{i}(G)$ for $i\in \{1,2,...,n\}$. A vertex with degree $k$ in $G$ will be referred as a $k$-vertex of $G$.  Throughout the paper, we always suppose that  $v_1$ is a vertex of $G$ with  $d_G(v_1)=d_1(G)$, and let  $n_k(G)$ define  the number of $k$-vertices of   $V(G)\setminus \{v_1\}$, namely, $$\text{$n_k(G)=|\{u:\, d_{G}(u)=k$ and $u\in V(G)\setminus \{v_1\}\}|$}.$$
	
Usually, $P_n$, $C_n$, $K_n$ and $K_{s,n-s}$  denote  the path,  the cycle, the complete graph and complete bipartite graph    with $n$ vertices, respectively. Especially, $K_{1,n-1}$ is a star with $n$ vertices.  For two disjoint graphs $G$ and $H$, we denote their disjoint union by $G\cup H$. And $G\vee H$ denotes the graph obtained from $G\cup H$ by adding all edges $vu$ with $v\in V(G)$ and $u\in V(H)$. Particularly, the disjoint union of $k$ copies of $G$ is denoted by $kG$.

Let $A(G)$ and $D(G)$ be the adjacency matrix and the diagonal degree matrix of $G$, respectively. The {\bf Laplacian    matrix} and {\bf signless  Laplacian    matrix} of $G$ are  defined as $L(G) = D(G)-A(G)$ and $Q(G) = D(G)+A(G)$, respectively.  For graph $G$ with $n$ vertices, let  $\kappa_{i}(G)$ and  $\lambda_{i}(G)$, respectively,   be the $i$-th largest signless Laplacian and  Laplacian eigenvalue of   $G$, where $1\le i\le n$. Let  $S_{Q}(G)$ and $S_{L}(G)$  denote the set of  eigenvalues of
$Q(G)$ and $L(G)$, respectively.
	It is easy to see that $Q(G)$ and $L(G)$ are  positive semidefinite and therefore their  eigenvalues can be, respectively,  arranged as \cite{ME1}
	$$\kappa_{1}(G)\geq \kappa_{2}(G)\geq\cdots\geq \kappa_{n}(G)\ge  0, \,\,\text{and}\,\,\lambda_{1}(G)\geq\lambda_{2}(G)\geq\cdots\geq\lambda_{n}(G)= 0$$ where  $\lambda_{n-1}(G)> 0$ if and only if $G$ is connected. If there is no risk of  confusion, then   $d_{G}(u)$, $N_{G}(u)$, $d_{i}(G)$, $\lambda_{i}(G)$ and $\kappa_i(G)$ are sometimes simplified   as $d(u)$, $N(u)$, $d_{i}$, $\lambda_{i}$ and $\kappa_i$, respectively.
	
	Two graphs are said to be {\bf $L$-cospectral} (resp., {\bf $Q$-cospectral})  if
	they have the same Laplacian   spectrum (resp., signless Laplacian   spectrum).
	Two graphs are said to be {\bf co-degree}  if they have the same degree sequences. A graph is called {\bf $DLS$} if there is no other non-isomorphic graph being $L$-cospectral with it.

 Identifying graphs that are, or are not, determined by the spectrum of an associated matrix is one of the oldest and the most extensively studied problems in the entire spectral theory \cite{Haemer2003,Haemer2009,Wang2017}. To the best of our knowledge, it dates back to 1950s when  G\"{u}nthard and Primas  considered it in the context  of the H\"{u}ckel's theory in Chemistry~\cite{GP}. In this line, Liu et  al. \cite{Liu2018} put forward the problem for: Which    cospectral graphs  have same degree sequences? In \cite{ZGL2017}, Liu et al. showed two sufficient conditions for two $Q$-cospectral graphs to be co-degree by setting ranges for signless Laplacian eigenvalues $\kappa_1$, $\kappa_2$ and $\kappa_n$, which is the main motivation of this paper, that is

\begin{theorem}\label{10t} {\em \cite{ZGL2017}} Let $G$ be a graph with $n\geq 12$ vertices and  $1\leq \kappa_{n}(G)\leq \kappa_{2}(G)\leq 5<n<\kappa_{1}(G)$.   If $S_{Q}(G)=S_{Q}(H)$, then $G$ and $H$ are  co-degree.
\end{theorem}
 \begin{theorem}\label{11t}{\em \cite{ZGL2017}} Let $G$ be a graph with  $n\geq 10$ vertices and  $0< \kappa_{n}(G)\leq \kappa_{2}(G)\leq 4<n<\kappa_{1}(G)$.   If $S_{Q}(G)=S_{Q}(H)$, then $G$ and $H$ are   co-degree.
\end{theorem}

It is natural to think what about Laplacian spectrum. Let $W_3$ and $W_5$ be the two graphs as shown in Fig. \ref{fig2}. We mainly show Theorem \ref{12t} and Theorem \ref{13t} in this paper as follows:

 \begin{theorem}\label{12t} Let $G$ be a connected  graph with $n\ge 18$ vertices and $\lambda_2(G)<5<n-1< \lambda_1(G)$. If $\lambda_{1}(G) \notin \{\lambda_{1}(W_3),\lambda_{1}(W_5)\}$ and $S_{L}(H)=S_{L}(G)$,  then $H$ and $G$ are  co-degree.
	\end{theorem}

 \begin{theorem}\label{13t} Let $G$ be a connected  graph with $n\ge 16$ vertices and $\lambda_2(G)\le 4.7<n-2< \lambda_1(G)$. If $S_{L}(H)=S_{L}(G)$,  then $H$ and $G$ are  co-degree.
	\end{theorem}

Let $G$ be a graph such that  every $L$-cospectral graph of  $G$ is also co-degree with  $G$. If the  degree sequence of $G$ is $\pi$ and every graph $H$ with $\pi$ as its degree sequence can be easily characterized, and we can easily prove that $G$ and $H$ must be isomorphic, then $G$ is  $DLS$. From this fact and Theorem \ref{12t}, we can easily  deduce the following conclusion, which is the  unique theorem result of \cite{LiuXG-DM}.

\begin{corollary}\label{c1.3}{\em \cite{LiuXG-DM}} Every multi-fan graph $K_1\vee (P_{l_1}\cup P_{l_1}\cup\cdots \cup P_{l_t})$ $(t\ge 1)$ is $DLS$.
\end{corollary}

What's more, we can also deduce the following conclusion:
\begin{corollary}\label{c1.4} Let $G:=K_1\vee (P_{l_1}\cup P_{l_1}\cup\cdots \cup P_{l_t}\cup C_{s_1}\cup C_{s_2}\cup\cdots \cup C_{s_k})$ $(t\ge 1, k\ge 1)$ with $n\ge 18$ vertices. If  each $s_i$ $(i=1,2,\ldots, k)$ is odd, then $G$ is $DLS$.
\end{corollary}

\section{Initial results}

 For a graph $G$ and vertex set  $X\subseteq V(G)$, we   write $G[X]$ to denote the subgraph induced by $X$ and  denote by $G-X$ the subgraph of $G$ induced by $V(G)\setminus X$. If $ X$ is a singleton vertex  $\{v\}$, then we simply write $G-v$. What's more, for $u,v \in V(G)$, we   use  $G-uv$ ($G+uv$) to denote the graph  obtained from $G$ by deleting (adding) the edge $uv$. And we let $a^{(k)}$ denote $k$ copies of number $a$.

\begin{lemma}\label{22l} {\em\cite{Guo2007}}
Let $U(G;tP_{k})$ be the graph obtained from a connected graph  $G$ by attaching the pendant  vertex $w_{i1}$ of $t$ paths $P_k$ to one non-isolated vertex $v$ of $G$, where the $i$-th path is $P^{(i)}_{k}=w_{i1}w_{i2}\cdots w_{ik}$ and $i\in \{1,2,\ldots,t\}$.   Let $\widehat{U}(G;tP_{k})$ be a  graph obtained from   $U(G;tP_{k})$ by adding any  $s$ edges among $\{w_{1k}, w_{2k},\ldots,w_{tk}\}$, where  $1\leq s\leq \frac{t(t-1)}{2}$.  Then,  $\lambda_1(U(G;tP_{k}))=\lambda_1(\widehat{U}(G;tP_{k}))$.
\end{lemma}

\begin{lemma}\label{25l}{\rm\cite{Das1}} Let $G$ be a graph with $N_{G}(u_{1})=N_{G}(u_{2})=\cdots=N_{G}(u_{k})=\{w_{1},w_{2},\ldots,w_{p}\}$. If $G^{*}$ is the graph obtained from $G$ by adding any $s$ $(1\leq s\leq \frac{k(k-1)}{2})$ edges among $\{u_{1},u_{2},\ldots,u_{k}\}$, then the eigenvalues of $L(G^{*})$ are as follows: those eigenvalues of $L(G)$ which are equal to $p$ $(k-1$ in number$)$ are incremented by $\lambda_{i}(G^{*}[X])$, $i=1,2,\ldots,k-1 $, and the remaining eigenvalues are the same, where $X=\{u_{1},u_{2},\ldots,u_{k}\}$.
\end{lemma}

\begin{lemma}\label{l-complement}{\rm\cite{Kelmans}} Let $\bar{G}$ be the complement of a graph $G$. If $S_L(G)=\{\lambda_1,\lambda_2,\cdots,\lambda_{n-1},0\}$, then $S_L(\bar{G})=\{n-\lambda_1,n-\lambda_2,\cdots,n-\lambda_{n-1},0\}$.
\end{lemma}

From Lemma \ref{l-complement}, we can easily deduce the Laplacian spectrum of $K_1\vee G$.

\begin{corollary}\label{c-join} Let $\lambda_1\ge \lambda_2\ge \cdots\ge \lambda_n(=0)$ be the Laplacian spectrum of  graph $G$. Then $n+1\ge 1+\lambda_1\ge 1+\lambda_2\ge \cdots\ge 1+\lambda_{n-1}\ge \lambda_n(=0)$ is the Laplacian spectrum of graph $K_1\vee G$.
\end{corollary}

From Corollary \ref{c-join}, we can further obtain the following result.

\begin{corollary}\label{c-multifan} Let $G:=K_1\vee (P_{l_1}\cup P_{l_1}\cup\cdots \cup P_{l_t}\cup C_{s_1}\cup C_{s_2}\cdots \cup C_{s_k})$ $(t\ge 1, k\ge 1)$ with $n$ vertices, then $S_L(G)$ is equal to\\
$$\left\{n,1^{(t+k-1)},0, 3-2\cos\left(\frac{i}{l_j}\pi\right): 1\le i\le l_j-1, 1\le j\le t, 3-2\cos\left(\frac{2i\pi}{s_{j}}\right): 1\leq i\leq s_{j}-1,  1\leq j\leq k\right\}.$$\\ Especially, if $G:=K_1\vee (P_{l_1}\cup P_{l_1}\cup\cdots \cup P_{l_t})$ $(t\ge 1)$ with $n$ vertices, then $S_L(G)$ is equal to\\
$$\left\{n,1^{(t-1)},0, 3-2\cos\left(\frac{i}{l_j}\pi\right): 1\le i\le l_j-1, 1\le j\le t\right\}.$$
\end{corollary}

\vspace*{0.2cm}

\begin{lemma}\label{L1-n}{\rm\cite{Liu-L1-n}} Let $G$ be a graph with $n$ vertices. The Laplacian spectrum of $G$ contains $n$ if and only if $G$ has a spanning subgraph that is a complete bipartite graph.
\end{lemma}

Let $uv$ be an edge of $G$  and $v$ be a vertex of $G$.  We denote by  $m(v)$  the average degree of the vertices in $N(v)$,  i.e., $m(v)=\underset{w\in
N(v)}{\sum}d(w)/d(v)$. We set
$$\Psi(uv)=\frac{d(u)(d(u)+m(u))+d(v)(d(v)+m(v))}{d(u)+d(v)}.$$

	\begin{lemma}\label{13l}{\rm\cite{ME1,Das1,pan2002sharp}} If $G$ is a connected graph with at least one edge, then \begin{align*}d_{1}(G)+1\leq  \lambda_1(G) \le \kappa_1(G)\leq \max\{\Psi(uv): uv\in E(G)\} \le     d_1(G)+d_2(G),\end{align*} where the first  equality holds   if and only if $d_1(G)=|V(G)|-1$,  the second  equality holds if and only if   $G$ is bipartite, and the third equality holds if and only if  $G$ is regular or  bipartite
semiregular.
	\end{lemma}

	\begin{lemma}\label{12l}{\rm\cite{pan2000second}}
		Let $G$ be a connected graph  with $ n\ge 3 $ vertices, then $\lambda_{2}(G)\ge d_{2}(G)$.
	\end{lemma}	

	 \begin{lemma}\label{l-d2=2}{\em \cite{Liu2018}} Let $G$ be a connected  graph with $n \ge 8$ vertices and $d_2(G)=2$. If $S_{L}(G)=S_{L}(H)$  then $G$ and $H$ are  co-degree.
	\end{lemma}

 The following interacting theorem for $G$ and $G-e$ is famous and useful. 	
\begin{lemma}\label{11l}{\rm\cite{van1995hamilton} }If $G$ is a graph with $ n\ge 3 $ vertices, then $0= \lambda_{n} (G-e)= \lambda_{n} (G) \le \lambda_{n-1}(G-e) \le \lambda_{n-1}(G) \le \dots \le \lambda_{1}(G-e) \le \lambda_{1}(G),$ and
$0\le  \kappa_{n} (G-e)\le \kappa_{n} (G) \le \kappa_{n-1}(G-e) \le \kappa_{n-1}(G) \le \dots \le \kappa_{1}(G-e) \le \kappa_{1}(G).$
\end{lemma}

\begin{figure}[H]
\vspace*{-11cm}\begin{center} \includegraphics[scale=1.0]{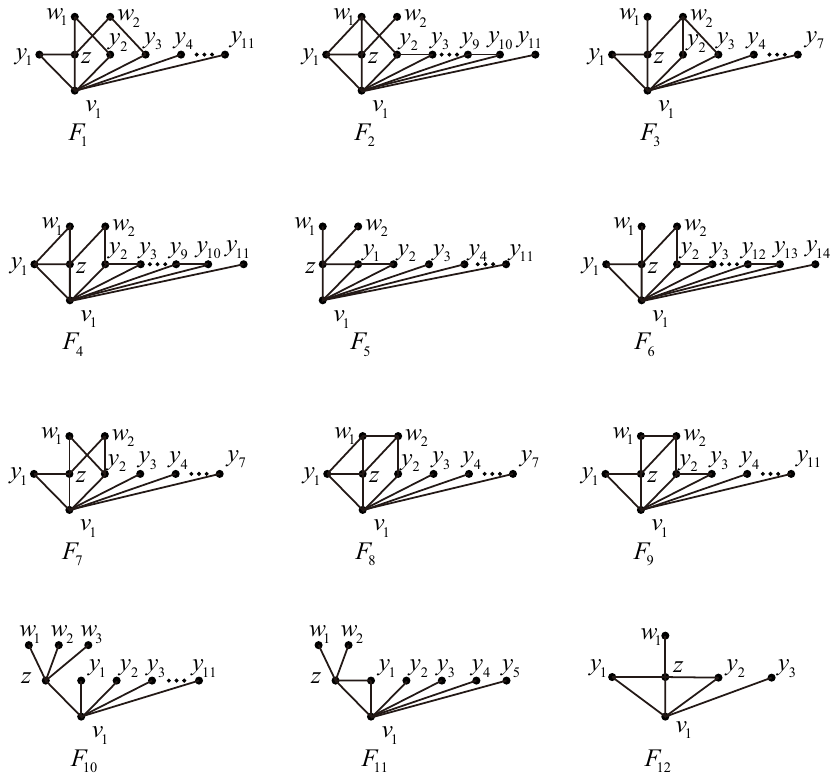} \par\vspace*{-3cm}
\caption{The   graphs    $F_{1}$, $F_2$, $\ldots$, $F_{12}$.} \label{fig3}\end{center}
     \end{figure}

	For graphs $G$ and $F$, we write $F\subset G$ if $F$ is a subgraph of $G$, and we call $F$ is {\bf forbidden} to $G$ if $G$ cannot contain $F$ as its subgraph. Let $F_1, F_2,\ldots, F_{12}$ be the $12$  graphs as shown in Figure \ref{fig3}. By direct calculations, $\lambda_2(F_i)\ge 5$, $i=1, 2,\ldots, 9$, and $\lambda_2(F_i)> 4.7$, $i=10, 11, 12$. Thus, if  $G$ is a graph  with $\lambda_2(G)<5$, then  \begin{align}\label{e3.1}\text{ each graph of  $\{$$F_i$: $1\le i\le 9$$\}$ is forbidden to $G$.}\end{align}
If  $G$ is a graph  with $\lambda_2(G)\le 4.7$, then  \begin{align}\label{e3.2}\text{ each graph of  $\{$$F_i$: $10\le i\le 12$$\}$ is forbidden to $G$.}\end{align}

\begin{figure}[H]
\vspace*{-11cm}\begin{center} \includegraphics[scale=1.0]{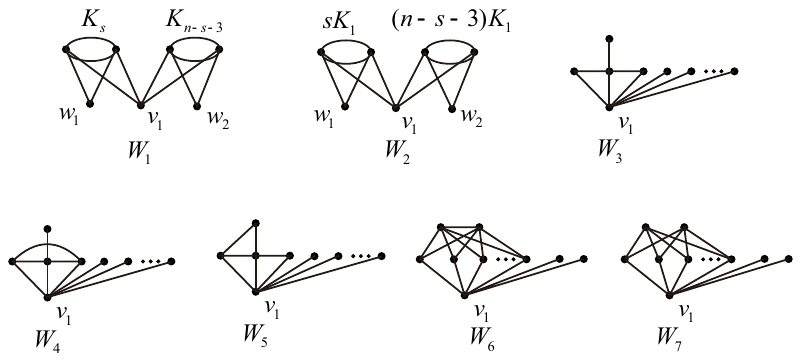} \par\vspace*{-11cm}
\caption{The   graphs    $W_{1}$, $W_2$, $\ldots$, $W_{7}$.} \label{fig2}\end{center}
     \end{figure}

Moreover, let $W_1, W_2,\ldots, W_{7}$ be the $7$ graphs as shown in Figure \ref{fig2}. We show an upper bound and a lower bound of the first largest Laplacian eigenvalue of graphs $W_1$ and $W_2$, that is:

\begin{lemma}\label{UpperBoundW1}
Let $W_1$ and $W_2$ be the graphs with $n\ge 5$ vertices in Fig \ref{fig2}, where  $1\le s\le n-4$. Then $n-2< \lambda_1(W_1)=\lambda_1(W_2)<n-1$.
\end{lemma}
\noindent {\bf Proof.} Let $v_1,w_1,w_2$ be the three vertices labeled in $W_2$ in Fig \ref{fig2}, where $d_{W_2}(v_1)=n-3$, $d_{W_2}(w_1)=s$ and $d_{W_2}(w_2)=n-3-s$. The remaining vertices all have degree $2$. By Lemmas \ref{25l} and \ref{13l}, $\lambda_1(W_1)=\lambda_1(W_2)>n-2$. Let $E_1$ and $E_2$ be the sets of edges incident with $w_1$ and $w_2$ in graph $W_2$, respectively. For the edges in $E(W_2)\setminus(E_1\cup E_2)$ (all incident with $v_1$), we let $E_3$ and $E_4$ be the sets of edges  adjacent to  the edges in $E_1$ and $E_2$, respectively. Then $E(W_2)=E_1\cup E_2 \cup E_3\cup E_4$.

If $uv\in E_1$, then
$$\Psi(uv)=\frac{s^2+2s+4+(s+n-3)}{s+2}=\frac{s^2+3s+n+1}{s+2}<n-1.$$

If $uv\in E_2$, then
$$\Psi(uv)=\frac{(n-3-s)^2+2(n-3-s)+4+(n-3+n-3-s)}{n-3-s+2}=\frac{(n-1-s)^2+s}{n-1-s}<n-1.$$

If $uv\in E_3$, then
$$\Psi(uv)=\frac{(n-3)^2+2(n-3)+4+(s+n-3)}{n-3+2}=\frac{(n-3)^2+3(n-3)+4+s}{n-1}<n-1.$$

If $uv\in E_4$, then
$$\Psi(uv)=\frac{(n-3)^2+2(n-3)+4+(n-3-s+n-3)}{n-3+2}=\frac{(n-3)^2+4(n-3)+4-s}{n-1}<n-1.$$
\vspace*{0.2cm}
Therefore, by Lemmas \ref{25l} and \ref{13l}, $\lambda_1(W_1)=\lambda_1(W_2)<n-1$. \qed
\vspace*{0.2cm}

   For a connected graph $G$, denote by $\tau(G)$ the number of triangles of  $G$. The following lemma is useful.
\begin{lemma}\label{14l}{\rm\cite{liu2012some}}
			If  $G$ is  a   graph with $n$ vertices and  $m$ edges,  then  \begin{align*}\sum_{i=1}^{n}\lambda_{i}=\sum_{i=1}^{n}d_{i}=2m, \hspace{5pt}\sum_{i=1}^{n}\lambda^{2}_{i}=2m+\sum_{i=1}^{n}d^{2}_{i},\,\text{and}\,\,
			&\sum_{i=1}^{n}\lambda^{3}_{i}=3\sum_{i=1}^{n}d^{2}_{i}+\sum_{i=1}^{n}d^{3}_{i}-6\tau(G).\end{align*}
		\end{lemma}

  \section{Proof of Theorem \ref{12t}}
  Let $\diag\{b_1, b_2, ... , b_n\}$ be a diagonal matrix with the given vectors lying on the main diagonal. And we denote by $\theta_{i}(M)$  the $i$-th largest eigenvalue of the square matrix $M$ with order $n$, where $1\le i\le n$.

	\begin{lemma}\label{l5.1}{\rm \cite{Liu2017}}
		Let $G$ be a graph with $\left\{w_{1},w_{2},\dots ,w_{s}\right\} $=W $ \subseteq V(G) $. If $ H\cong G[W] $, then for
		any $ i\in\left\{1,2,\dots,s\right\} $ and $ 0\le p_{i} \le d_{G}(w_i) $, we
		have $ \lambda_{i}(G)\ge \theta_{i}(diag\left\{p_{1},p_{2},\dots\,p_{s}\right\}-A(H)) $.
	\end{lemma}

	\begin{lemma}{\rm\cite{Ye2024two}}\label{e7}
		Let $G$ be a connected graph with $n$ vertices and $ d_2(G)\le 4 $. If either $d_n(G)=1<11\le d_1(G)$ or $2\le d_n(G)<8\le d_1(G)$, then $\lambda_1(G)\le \kappa_1(G)\le d_1(G)+3$.
	\end{lemma}

\begin{figure}[H]
\vspace*{-12cm}\begin{center} \includegraphics[scale=1.0]{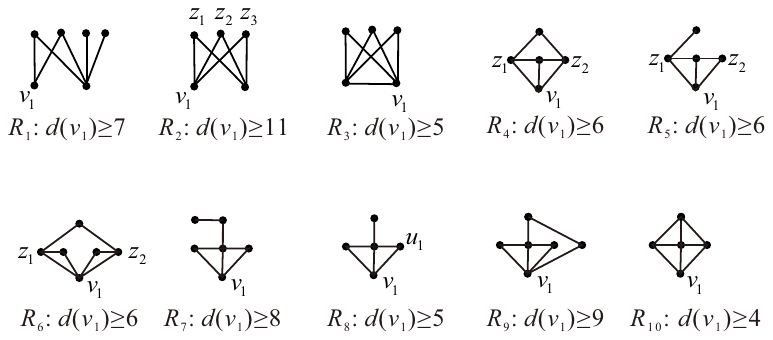} \par\vspace*{-11.5cm}
\caption{The   graphs    $R_{1}$, $R_2$, $\ldots$, $R_{10}$.} \label{fig1}\end{center}
     \end{figure}

  Next, we will use vertex-weighted graphs illustrated in Figure~\ref{fig1}. They are `vertex-weighted' in the sense that some of their vertices are accompanied with weights. Precisely, vertices $z_1$, $z_2$ and $z_3$  have weight 4, vertex $u_1$ has weight $3$, and weight $d(v_1)$ of  $v_1$ is as in the figure. These graphs can be seen as subgraphs of larger graphs, say $R'_1, R'_2, \ldots, R'_{10}$, in which degrees of the previously mentioned vertices are equal to given weights. We don't need to know the exact distribution of  other edges, but we care about the vertex degrees, and that is why they are interpreted as weights. More succinctly, we use the same  notation  for vertex weights and vertex degrees. Therefore, for $1\leq i\leq 10$, the matrix representation of $R_i$, say $L^{*}(R_i)$, is the principal submatrix of the Laplacian matrix $L(R'_i)$.

  \begin{remark} In what follows, we consider graphs $G$  that contain some of graphs $R'_i, 1\leq i\leq 10$, as subgraphs, and in each situation we estimate $\lambda_2(G)$ as demonstrated in this remark. Take, for example, that a graph $G$ contains $R'_8$ as a subgraph $($Notice that $u_1$ has weight $3$ in $R_8$ in Fig \ref{fig1}, which means that there exists some vertex, say $z$, such that $z\in V(R'_8)\setminus V(R_8)$ and $zu_1\in E(R'_8)$$)$. Then $L^{*}(R_8)$ is a principal submatrix of $L(R'_8)$. By Lemmas \ref{11l} and \ref{l5.1}, we have  $\lambda_2(G)\ge \theta_2(L^*(R_8))>5$ when $d_{R'_8}(v_1)\ge 5$, where 			\begin{center}
				$L^*(R_8) =	\begin{pmatrix}
					&d_{R'_8}(v_1)&-1&-1&-1 &0\\
					& -1& 4 &-1 &-1 &-1\\
				    &-1 &-1 & 3 & 0 & 0\\
					&-1 &-1 & 0 & 2 & 0\\
					&0 &-1 & 0 & 0 & 1
				\end{pmatrix} $.
			\end{center}
  	\end{remark}
		
\begin{lemma}\label{l5.3}
Let $G$ be a connected graph with $n\ge15$ vertices such that $\lambda_{2}(G)<5$ and $d_{1}(G)\ge n-4$. Then, $d_{2}(G)\le 4$ and  any two $4$-vertices (if exist) are not adjacent. Moreover, if $w\notin \{v_1\}\cup  N_G(v_1)$, then $d_{G}(w)\le 3$ and  $w$ is adjacent to at most two  $4$-vertices.
\end{lemma}
\noindent {\bf Proof.}
		By Lemma \ref{12l}, we have $d_2(G)\le 4$ since $\lambda_{2}(G)<5$.  Let $w$ be a vertex with $w\notin \{v_1\}\cup  N_G(v_1)$ (if exists).
 If $d_G(w)=4$,  then by Lemmas \ref{11l} and \ref{l5.1} we have $\lambda_{2}(G)\ge \theta_2({L}^*(R_1))\ge 5 $, a contradiction. Thus $d_G(w)\le 3 $. This implies that every 4-vertex is adjacent to $v_1$.   Assume that there exists two 4-vertices being adjacent, then the two 4-vertices and $v_1$ induces a triangle. One can easily check that $\lambda_2(G)\ge 5$ when $d_G(v_1)\ge 4$ by Lemma \ref{l5.1}, a contradiction. So we can conclude that any two $4$-vertices (if exist) are not adjacent.

  Next, we show  that $w$ is adjacent to at most two 4-vertices. Otherwise, suppose that $d_G(w)=3$ and $w$ is adjacent to three 4-vertices. Then, the three 4-vertices induces an independent set of $G$, and then by Lemmas \ref{11l} and \ref{l5.1} we have $\lambda_{2}(G)\ge \theta_{2}(L^*(R_2))\ge 5 $, a contradiction. \qed
\begin{lemma}\label{3.8}
	Let $G$ be a connected graph with $n\ge 15$ vertices, $d_2(G)=4$ and   $\lambda_{2}(G)<5<n-1<\lambda_1(G)$. Then $n-3\le d_1(G)\le n-2$ and  $n_4(G)=1$. Moreover, $zv_1\in E(G)$ and $1\le |N_G(z)\cap N_G(v_1)|\le 2$, where $z$ is the unique $4$-vertex of $G$.
	\end{lemma}
\noindent{\bf Proof.} Since  $d_2(G)=4$ and $n-1<\lambda_1(G)\le d_1(G)+4$ (by Lemma \ref{13l}), we have $d_1(G)\ge n-4\ge 11$. In view of Lemma \ref{e7},  it follows that   $d_1(G)\ge n-3$.  Let  $z$ be a 4-vertex of $G$. By Lemma \ref{l5.3}, we have  $zv_1\in E(G)$. Note that $d_G(v_1)\ge n-3$, so $\left| N_G(z)\cap N_G(v_1) \right|\ge 1$. If $\left| N_G(z)\cap N_G(v_1) \right|=3$, then by Lemmas \ref{11l} and \ref{l5.1}, we have $\lambda_2(G)\ge \theta_{2}(L^{*}(R_3))\ge 5 $, a contradiction. Thus, $1\le |N_G(z)\cap N_G(v_1)|\le 2$ and $n-3\le d_1(G)\le n-2$.

Next, we show that $n_4(G)=1$. By contradiction, assume that $n_4(G)\ge 2$ and $z_1, z_2$ are two 4-vertices of $G$. By Lemma \ref{l5.3}, $z_1z_2\not\in E(G)$ and $\{z_1, z_2\}\subset N_G(v_1)$. Since $1\le | N_G(z_1)\cap N_G(v_1)|\le 2$, we can suppose $w_1\in N_G(z_1)\setminus (N_G(v_1)\cup\{v_1\})$. If there exists $u\in N_G(v_1)$ such that $u\in N_G(z_1)\cap N_G(z_2)$, then either  $\lambda_2(G)\ge\theta_2(L^{*}(R_4))\ge 5$ ($w_1z_2\in E(G)$) or $\lambda_2(G)\ge\theta_2(L^{*}(R_5))\ge 5$ ($w_1z_2\notin E(G)$), a contradiction. Thus, $N_G(z_1)\cap N_G(z_2)\cap N_G(v_1)=\emptyset$.

Now, note that $1\le |N_G(z_2)\cap N_G(v_1)|\le 2$, if $w_1z_2\in E(G)$, then $\lambda_2(G)\ge\theta_2(L^{*}(R_6))\ge 5$, a contradiction. So we can suppose $w_2\in N_G(z_2)\setminus (N_G(v_1)\cup\{v_1\})$ and $w_2\neq w_1$. And hence $d_1(G)=n-3$. If $n_4(G)\ge3$, we assume that $z_3$ is an another $4$-vertex. Then from above statements, we have either $w_1z_3\in E(G)$ or $w_2z_3\in E(G)$. And so $\lambda_2(G)\ge\theta_2(L^{*}(R_6))\ge 5$, a contradiction. Thus $n_4(G)=2$. Now we would like to show that $N_G(w_i)\setminus\{z_i\}\subset N_G(z_i)$, $i\in\{1, 2\}$. Otherwise, $\lambda_2(G)\ge\theta_2(L^{*}(R_7))\ge 5$, a contradiction. Note that $N_G(z_1)\cap N_G(z_2)=\{v_1\}$. It is easy to see that $w_1$ and $w_2$ are in different components of $G-v_1$ (Otherwise $\lambda_2(G)\ge\theta_2(L^{*}(R_8))\ge 5$). So $G$ is a subgraph of $W_{1}$ in Fig \ref{fig2}. By Lemmas \ref{11l} and \ref{UpperBoundW1}, $\lambda_1(G)\le \lambda_1(W_{1})=\lambda_1(W_{2})< n-1$, a contradiction. This confirms $n_4(G)=1$. \qed

		\par\medskip
\noindent{\bf Proof of Theorem \ref{12t}.}  Since $G$ is connected, $\lambda_{n-1}(G)> 0$. Then $\lambda_{n-1}(H)=\lambda_{n-1}(G)> 0$, so $H$ is also connected.  For short,  next we write $n_i(G)$, $d_j(G)$, $n_i(H)$ and $d_j(H)$ as $n_i$, $d_j$, $n_i^*$ and $d_{j}^*$, respectively, where $i\in \{1,2,3,4\}$ and $j\in \{1,2,\ldots, n\}$.  Without loss of generality, we suppose $d_1\le d^*_1$.
	By   Lemmas \ref{12l} and \ref{13l}, we have $d_2\le 4$, $d_{2}^*\le4$, $d_1\ge n-4$ and $d_{1}^*\ge n-4$. Since $n-4>11$, it follows from Lemma \ref{e7} that   $d_1\ge n-3$ and $d_{1}^*\ge n-3$. Since $K_{1,n-1}$ is determined by its Laplacian spectrum (see \cite{Liu2010}), taking  Lemma \ref{l-d2=2} into consideration, we may suppose that   $d_2\ge 3$ and $d_{2}^*\ge 3$.

From Lemma \ref{14l}, we have 	
 \begin{equation}\label{e51}\left\{
	\begin{array}{ll}&  2n_{3}^*+6n_{4}^*+d_{1}^*(d_{1}^*-3)=2n_3+6n_4+d_1(d_1-3),\\ &n_{2}^*-3n_{4}^*-d_{1}^*(d_{1}^*-4)=n_2-3n_4-d_1(d_1-4),\\
&2n_{1}^*+2n_{4}^*+d_{1}^*(d_{1}^*-5)=2n_1+2n_4+d_1(d_1-5).
	\end{array}
	\right.\end{equation}

 We claim that $d_1=d^*_1$. Conversely we assume that $d_1< d^*_1$. Then,   $ n-3\le d_1<d^*_1\le n-1$.
			Note that $d_{2}^*\ge 3$, so  $n-2\ge n^*_2=\big(n_{2}+3n^*_{4}+d^*_{1}(d^*_{1}-4)-d_{1}(d_{1}-4)-3n_{4}\big)$ by \eqref{e51}. By Lemma  \ref{3.8}, we know $n_4\le 1$. And note that  there are three possibilities for $d_1$ and $d^*_1$: If $d_{1}^*=n-1$ and $d_1=n-2$, then  $n-2\ge    d^*_{1}(d^*_{1}-4)-d_{1}(d_{1}-4)-3=2n-10$. And hence $n\le 8$, contrary with $n\ge 18$; If $d_{1}^*=n-1$ and $d_1=n-3$, then  $n-2\ge    d^*_{1}(d^*_{1}-4)-d_{1}(d_{1}-4)-3=4n-19$. And hence $n\le 5$, contrary with $n\ge 18$; If $d_{1}^*=n-2$ and $d_1=n-3$, then  $n-2\ge    d^*_{1}(d^*_{1}-4)-d_{1}(d_{1}-4)-3=2n-12$. And hence $n\le 10$, contrary with $n\ge 18$. This confirms  $d_1=d^*_1$.

Now combining $d_1=d^*_1$ with \eqref{e51}, we have
 \begin{equation}\label{e52}\left\{
	\begin{array}{ll}&  n_{3}^*+3n_{4}^*=n_3+3n_4,\\ &n_{2}^*-3n_{4}^*=n_2-3n_4,\\
&n_{1}^*+n_{4}^*=n_1+n_4.
	\end{array}
	\right.\end{equation}

 From \eqref{e52}, if $n_4=n^*_4$, then $n_i=n_{i}^*$, $i=1, 2 ,3$. The result already holds. Next, without loss of generality, we conversely assume that $n_4>n^*_4$.  By Lemma  \ref{3.8}, we know $n_{4}\le 1$. So $0=n^*_4<n_4=1$.  Combining this with \eqref{e52}, we have
    \begin{equation}\label{e53}\left\{
	\begin{array}{ll}& n_{3}^*=n_3+3,\\ &n_{2}^*=n_2-3,\\
&n_{1}^*=n_1+1.
	\end{array}
	\right.\end{equation}

   By Lemma \ref{3.8}, we obtain that $n-3\le d^*_1=d_1\le n-2$ as $n_4=1$.

\par\noindent{\bf Case 1. $d_1=d^*_1=n-2$.} Suppose that $z$ is  the unique 4-vertex of $G$, $d_G(v_1)=n-2$ and $w_1=V(G)\setminus(N_G(v_1)\cup\{v_1\})$. By Lemma \ref{l5.3}, we have $zv_1\in E(G)$.  Obviously, $zw_1\in E(G)$ (Otherwise $\lambda_2(G)\ge \theta_2(L^*(R_3))\ge 5$). Suppose that $N_G(z)=\{v_1, w_1, w_2, w_3\}$. By Lemma \ref{l5.3}, $1\le d_G(w_1)\le3$ and   $2\le d_G(w_i)\le3$, $i=2, 3$.

\par\noindent{\bf Subcase 1.1. $d_G(w_1)=1$.} $(i)$: If $d_G(w_2)=d_G(w_3)=2$, then by Lemma \ref{22l} or Lemma \ref{25l}, $\lambda_1(G)=\lambda_1(W_3)$, a contradiction. $(ii)$: If $d_G(w_2)=3>d_G(w_3)=2$ or $d_G(w_2)=2<d_G(w_3)=3$, then by Lemma \ref{l5.1}, $\lambda_2(G)\ge\theta_2(L^{*}(R_8))\ge 5$, a contradiction. $(iii)$: $d_G(w_2)=d_G(w_3)=3$. If $w_2w_3\notin E(G)$, then by Lemmas \ref{11l} and \ref{l5.1}, $\lambda_2(G)\ge\theta_2(L^{*}(R_8))\ge 5$, a contradiction. Otherwise, $w_2w_3\in E(G)$. But now by Lemmas \ref{22l} and  \ref{25l}, $\lambda_1(G)=\lambda_1(W_4)=\lambda_1(W_3)$, a contradiction.

\par\noindent{\bf Subcase 1.2. $d_G(w_1)=2$.} $(i)$: If $w_1w_2\in E(G)$ or $w_1w_3\in E(G)$, by symmetry we just consider $w_1w_2\in E(G)$. $(a)$: If $d_G(w_3)=3$, then $\lambda_2(G)\ge\lambda_2(G-w_1w_2)\ge\theta_2(L^{*}(R_8))\ge 5$, a contradiction. $(b)$: If $d_G(w_3)=2$, then by Lemma \ref{22l}, $\lambda_1(G)=\lambda_1(W_5)$, a contradiction. $(ii)$: If $w_1w_2\notin E(G)$ and $w_1w_3\notin E(G)$, then $\lambda_2(G)\ge\theta_2(L^{*}(R_9))\ge 5$, a contradiction.

\par\noindent{\bf Subcase 1.3. $d_G(w_1)=3$.} $(i)$: If $w_1w_2\notin E(G)$ and $w_1w_3\notin E(G)$, then $\lambda_2(G)\ge\theta_2(L^{*}(R_9))\ge 5$, a contradiction.
$(ii)$: If exactly one vertex of $\{w_2, w_3\}$ is adjacent to $w_1$, then  by Lemmas \ref{11l} and \ref{l5.1}, we have $\lambda_2(G)\ge\theta_2(L^{*}(R_9))\ge 5$, a contradiction. $(iii)$: If $w_1w_2\in E(G)$ and $w_1w_3\in E(G)$, then $\lambda_2(G)\ge\theta_2(L^{*}(R_{10}))\ge 5$, a contradiction.

 \par\noindent{\bf Case 2. $d_1=d^*_1=n-3$.}  Let $V(G)\setminus N_G(v_1)=\{w_1,w_2,v_1\}$.
If $w_1$ and $w_2$ are in different components of $G-v_1$, then $G$ is a subgraph of $W_{1}$. By Lemmas \ref{11l} and \ref{UpperBoundW1}, $\lambda_1(G)\le \lambda_1(W_{1})=\lambda_1(W_{2})< n-1$, a contradiction. Thus, $w_1$ and $w_2$ are  in the same  component, say $G_1$, of $G-v_1$.  Let  $G_2=(G-v_1)-V(G_1)$. By Lemma \ref{22l} , we have $\lambda_1(G)=\lambda_1(\widetilde{G})$, where $\widetilde{G}$ is obtained from $G$ by deleting all edges of $E(G_2)$. Note that $d_2=4$ and $n_4=1$. If $|V(G_2)|\ge 2$, by Lemmas \ref{25l} and \ref{11l}, then $\lambda_1(G)=\lambda_1(\widetilde{G})\le\lambda_1(W_{6})=\lambda_1(W_{7})<n-1$, a contradiction. Thus, $|V(G_2)|\le 1$ and $n-2\le |V(G_1)|\le n-1$.

Let $z$ be the unique 4-vertex of $G$. By Lemma \ref{l5.3}, we have $zv_1\in E(G)$. Moreover,  $1\le |N_G(z)\cap N_G(v_1)|\le 2$ by Lemma \ref{3.8}. So at least one of $\{w_1, w_2\}$ is adjacent to $z$, without loss of generality, we suppose  $w_1z\in E(G)$. Next we would like to show that $|N_G(z)\cap N_G(v_1)|=1$. By contradiction, we assume that $|N_G(z)\cap N_G(v_1)|=2$ and $N_G(z)\cap N_G(v_1)=\{y_1,y_2\}$. Now, obviously, $N_G(w_1)\subset \{z,y_1,y_2\}$, otherwise $\lambda_2(G)\ge\theta_2(L^{*}(R_7))\ge 5$, a contradiction. And  $N_G(y_i)\subset \{z,v_1,w_1\}$, $i=1, 2$, otherwise $\lambda_2(G)\ge\theta_2(L^{*}(R_8))\ge 5$, a contradiction. However, it means that $w_1$ and $w_2$ are in different components of $G-v_1$, a contradiction. So we can conclude that $|N_G(z)\cap N_G(v_1)|=1$, $zw_1\in E(G)$ and $zw_2\in E(G)$. For convenience, from now on we suppose $N_G(z)\cap N_G(v_1)=\{y_1\}$ and $N_G(v_1)=\{z, y_1, y_2,\ldots,y_{n-4} \}$.

\par\noindent{\bf Subcase 2.1. $w_1w_2\notin E(G)$.}
\par\noindent{\bf Subcase 2.1.1. $w_1y_1\in E(G)$ or $w_2y_1\in E(G)$.} Note that $n_4=1$ and $d_2=4$, so in this case, exactly one vertex of $\{w_1,w_2\}$ is adjacent to $y_1$. By symmetry, we may suppose that $w_1y_1\in E(G)$. Then $2\le d_G(w_1)\le 3$ and $1\le d_G(w_2)\le 3$.

$(i)$: $d_G(w_1)=3$. And we may suppose $w_1y_2\in E(G)$. $(a)$: If $d_G(w_2)=2$ and $w_2y_2\in E(G)$, then $\lambda_1(G)\le \lambda_1(W_7)<n-1$ since $d_2=4$ and $n_4=1$, a contradiction. $(b)$: If $d_G(w_2)=2$ and $w_2y_2\notin E(G)$, then $\lambda_2(G)\ge \lambda_2(F_1)\ge 5$, a contradiction. $(c)$: If $d_G(w_2)=3$,  then $\lambda_2(G)\ge \lambda_2(F_1)\ge 5$, a contradiction. $(d)$: Otherwise $d_G(w_2)=1$. Note that $n-2\le |V(G_1)|\le n-1$, $n_4=1$, $d_2=4$ and $n\ge 18$. So it is easy to see that $F_2\subset G$. And hence $\lambda_2(G)\ge \lambda_2(F_2)>5$, a contradiction.

$(ii)$: $d_G(w_1)=2$.  $(a)$: If $d_G(w_2)=3$,  then $\lambda_2(G)\ge \lambda_2(F_3)> 5$, a contradiction. $(b)$: If $d_G(w_2)=1$, then $\lambda_1(G)\le \lambda_1(W_7)<n-1$ since $d_2=4$ and $n_4=1$, a contradiction. $(c)$: Otherwise $d_G(w_2)=2$. Note that $n-2\le |V(G_1)|\le n-1$, $n_4=1$, $d_2=4$ and $n\ge 18$. So it is easy to see that $F_4\subset G$. And hence $\lambda_2(G)\ge \lambda_2(F_4)>5$, a contradiction.

\par\noindent{\bf Subcase 2.1.2. $w_1y_1\notin E(G)$ and $w_2y_1\notin E(G)$.} In this case, we claim that $d_G(y_1)=2$. Otherwise $d_G(y_1)=3$. Then $F_5\subset G$ and $\lambda_2(G)\ge \lambda_2(F_5)>5$, a contradiction. So $N_G(y_1)=\{z, v_1\}$.

$(i)$: $d_G(w_1)=d_G(w_2)=1$. Then $\lambda_1(G)\le \lambda_1(W_7)<n-1$, a contradiction.

$(ii)$: $\{d_G(w_1), d_G(w_2)\}=\{1,2\}$. By symmetry, we may suppose that $d_G(w_1)=1$ and $d_G(w_2)=2$. Note that $n-2\le |V(G_1)|\le n-1$, $n_4=1$, $d_2=4$ and $n\ge 18$. So it is easy to see that $F_6\subset G$. And hence $\lambda_2(G)\ge \lambda_2(F_6)>5$, a contradiction.

$(iii)$: $\{d_G(w_1), d_G(w_2)\}=\{1,3\}$. By symmetry, we may suppose that $d_G(w_1)=1$ and $d_G(w_2)=3$. Then it is easy to see that $F_3\subset G$. And hence $\lambda_2(G)\ge \lambda_2(F_3)>5$, a contradiction.

$(iv)$: $d_G(w_1)\ge2$ and $d_G(w_2)\ge 2$. If $|N_G(w_1)\cap N_G(w_2)|\ge 2$, then $F_7\subset G$ and $\lambda_2(G)\ge \lambda_2(F_7)>5$, a contradiction. Otherwise, $|N_G(w_1)\cap N_G(w_2)|=1$. Then $F_1\subset G$ and $\lambda_2(G)\ge \lambda_2(F_1)\ge 5$, a contradiction.

\par\noindent{\bf Subcase 2.2. $w_1w_2\in E(G)$.} In this case, $2\le d_G(w_i)\le 3$, $i=1,2$. As before, we will consider some subcases as following.

\par\noindent{\bf Subcase 2.2.1. $w_1y_1\in E(G)$ or $w_2y_1\in E(G)$.} Note that $n_4=1$ and $d_2=4$, so in this case, exactly one vertex of $\{w_1,w_2\}$ is adjacent to $y_1$. By symmetry, we may suppose that $w_1y_1\in E(G)$. Then $d_G(w_1)=3$ and $2\le d_G(w_2)\le 3$. If $d_G(w_2)=2$, then by Lemmas \ref{25l} and \ref{11l}, $\lambda_1(G)=\lambda_1(\widetilde{G})\le \lambda_1(\widetilde{G}+y_1w_2)\le \lambda_1(W_{6})=\lambda_1(W_{7})<n-1$, a contradiction. Otherwise, $d_G(w_2)=3$. Then $F_8\subset G$ and $\lambda_2(G)\ge \lambda_2(F_8)>5$, a contradiction.

\par\noindent{\bf Subcase 2.2.2. $w_1y_1\notin E(G)$ and $w_2y_1\notin E(G)$.}

$(i)$: $d_G(w_1)=d_G(w_2)=2$. If $d_G(y_1)=2$, then by Lemmas \ref{25l} and \ref{11l}, $\lambda_1(G)=\lambda_1(\widetilde{G})\le \lambda_1(W_{6})=\lambda_1(W_{7})<n-1$, a contradiction. Otherwise, $d_G(y_1)=3$. Then $F_5\subset G$ and $\lambda_2(G)\ge \lambda_2(F_5)>5$, a contradiction.

$(ii)$: $\{d_G(w_1), d_G(w_2)\}=\{2,3\}$. By symmetry, we may suppose that $d_G(w_1)=2$, $d_G(w_2)=3$ and $w_2y_2\in E(G)$ ($y_2\notin \{y_1,w_1,z\}$). Note that $F_5\not \subset G$, so $d_G(y_1)\not=3$, i.e. $d_G(y_1)=2$.  Then it is easy to see that $F_{9}\subset G$ since $n-2\le |V(G_1)|\le n-1$, $n_4=1$, $d_2=4$ and $n\ge 18$. And hence $\lambda_2(G)\ge \lambda_2(F_{9})\ge 5$, a contradiction.

$(iii)$: $d_G(w_1)=d_G(w_2)=3$. If $|N_G(w_1)\cap N_G(w_2)|\ge 2$, then $F_7\subset G$ and $\lambda_2(G)\ge \lambda_2(F_7)>5$, a contradiction. Otherwise, $|N_G(w_1)\cap N_G(w_2)|=1$. Then $F_1\subset G$ and $\lambda_2(G)\ge \lambda_2(F_1)\ge 5$, a contradiction.

In summary, we can conclude that if  graphs G and H satisfy the conditions of Theorem \ref{12t}, then G and H must be co-degree.  \qed

In view of the proof of Theorem \ref{12t}, we can easily obtain the following corollary:
 \begin{corollary}\label{c3.1} Let $G$ be a connected  graph with $n\ge 18$ vertices and $\lambda_2(G)<5<n-1< \lambda_1(G)$. If $d_{1}(G) \neq n-2$ and $S_{L}(H)=S_{L}(G)$,  then $H$ and $G$ are  co-degree.
	\end{corollary}
\section{Proofs of Theorem \ref{13t} and Corollaries \ref{c1.3} and \ref{c1.4}}


\begin{lemma}\label{l6.1}
Let $G$ be a connected graph with $n\ge 16$ vertices. If $\lambda_{2}(G)\le 4.7<n-2<\lambda_1(G)$, then $d_2(G)\le 3$ and $d_1(G)\ge n-4$.
	\end{lemma}
\noindent{\bf Proof.} Since $\lambda_{2}(G)\le 4.7$, by Lemma \ref{12l}, $d_2(G)\le 4$. Then $n-2<\lambda_1(G)\le d_1(G)+4$ by Lemma \ref{13l}, so $d_1(G)\ge n-5\ge11$. Now, by Lemma \ref{e7}, we have $n-2<\lambda_1(G)\le d_1(G)+3$, so $d_1(G)\ge n-4$.

Next we will show that $d_2(G)\le 3$. By contradiction, we assume that $d_2(G)=4$ and $z$ is a 4-vertex. By Lemma \ref{l5.3}, $zv_1\in E(G)$. If $|N_G(z)\cap N_G(v_1)|=0$, then $F_{10}\subset G$, so $\lambda_2(G)\ge \lambda_2(F_{10})>4.7$, a contradiction; If $|N_G(z)\cap N_G(v_1)|=1$, then $F_{11}\subset G$, so $\lambda_2(G)\ge \lambda_2(F_{11})>4.7$, a contradiction; If $|N_G(z)\cap N_G(v_1)|=2$, then $F_{12}\subset G$, so $\lambda_2(G)\ge \lambda_2(F_{12})>4.7$, a contradiction; If $|N_G(z)\cap N_G(v_1)|=3$, then by Lemmas \ref{11l} and \ref{l5.1}, we have $\lambda_2(G)\ge \theta_{2}(L^{*}(R_3))\ge 5>4.7$, a contradiction. Therefore, $n_4(G)=0$ and $d_2(G)\le 3$. \qed

\par\medskip
\noindent{\bf Proof of Theorem \ref{13t}.}
By Lemma \ref{l6.1}, $d_2(G)\le 3$, $d_2(H)\le 3$, $d_1(G)\ge n-4$  and $d_1(H)\ge n-4$.
In the following, we simplify $n_i(G)$, $n_i(H)$, $d_1(G)$ and $d_1(H)$ as $n_i$, $n_i^*$, $d_1$ and $d_{1}^*$, respectively, where $i\in \{1,2,3\}$. Since $K_{1,n-1}$ is determined by its Laplacian spectrum (see \cite{Liu2010}), taking  Lemma \ref{l-d2=2} into consideration, we may suppose that   $d_2=d_{2}^*=3$. And without loss of generality, we suppose  $d_1\le d_1^*$. By Lemma \ref{14l}, we have
\begin{align}\label{e6.1}\left\{
\begin{array}{lll}
n^{*}_{1}+n^{*}_{2}+n^{*}_{3}=n_{1}+n_{2}+n_{3},\\
n^{*}_{1}+2n^{*}_{2}+3n^{*}_{3}+d^{*}_{1}=n_{1}+2n_{2}+3n_{3}+d_{1},\\
n^{*}_{1}+4n^{*}_{2}+9n^{*}_{3}+(d^{*}_{1})^2=n_{1}+4n_{2}+9n_{3}+(d_{1})^{2}.
\end{array}
\right.\end{align}
From \eqref{e6.1}, it is easy to get that
\begin{align}\label{e6.2}
n^*_2=n_{2}+d^*_{1}(d^*_{1}-4)-d_{1}(d_{1}-4).
\end{align}

We claim that $d_1=d^*_1$. Conversely we assume that $d_1< d^*_1$. Then,   $ n-4\le d_1<d^*_1\le n-1$.
			Since $d_{2}^*=3$, we get $n-2\ge n^*_2\ge d^*_{1}(d^*_{1}-4)-d_{1}(d_{1}-4)$ by \eqref{e6.2}. And note that  there are six possibilities for $d_1$ and $d^*_1$: If $d_{1}^*=n-1$ and $d_1=n-2$, then  $n-2\ge  d^*_{1}(d^*_{1}-4)-d_{1}(d_{1}-4)=2n-7$. And hence $n\le 5$, contrary with $n\ge 16$; If $d_{1}^*=n-1$ and $d_1=n-3$, then  $n-2\ge    d^*_{1}(d^*_{1}-4)-d_{1}(d_{1}-4)=4n-16$. And hence $n\le 4$, contrary with $n\ge 16$; If $d_{1}^*=n-1$ and $d_1=n-4$, then  $n-2\ge d^*_{1}(d^*_{1}-4)-d_{1}(d_{1}-4)=6n-27$. And hence $n\le 5$, contrary with $n\ge 16$; If $d_{1}^*=n-2$ and $d_1=n-3$, then  $n-2\ge    d^*_{1}(d^*_{1}-4)-d_{1}(d_{1}-4)=2n-9$. And hence $n\le 7$, contrary with $n\ge 16$; If $d_{1}^*=n-2$ and $d_1=n-4$, then  $n-2\ge    d^*_{1}(d^*_{1}-4)-d_{1}(d_{1}-4)=4n-20$. And hence $n\le 6$, contrary with $n\ge 16$; If $d_{1}^*=n-3$ and $d_1=n-4$, then  $n-2\ge    d^*_{1}(d^*_{1}-4)-d_{1}(d_{1}-4)=2n-11$. And hence $n\le 9$, contrary with $n\ge 16$. This confirms  $d_1=d^*_1$.

Now combining $d_1=d^*_1$ with \eqref{e6.2}, we get $n_2^*=n_2$. Then by \eqref{e6.1}, it is easy to see that $n_1^*=n_1$ and $n_3^*=n_3$. Therefore, $G$ and $H$ are co-degree. \qed
\par\medskip
\noindent{\bf Proof of Corollary \ref{c1.3}.} For convenience, let $G:=K_1\vee (P_{l_1}\cup P_{l_1}\cup\cdots \cup P_{l_t})=K_1\vee (sP_1\cup P_{l_{s+1}}\cup P_{l_{s+2}}\cup\cdots \cup P_{l_t})$ ($2\le l_{s+1}\le l_{s+2} \le\cdots \le l_t$), $n:=1+\sum_{i=1}^{t}l_{i}$ and $H$ be $L$-cospectral with $G$. By Corollary \ref{c-multifan} and Lemma \ref{L1-n}, $\lambda_1(G)=n\notin\{\lambda_1(W_3),\lambda_1(W_5)\}$  and $\lambda_2(G)=3+2\cos\pi/l_{t}<5$. When $n\le 17$, it is easy to check that $H=G$ by Lemmas \ref{L1-n}, \ref{13l}  and \ref{12l}. So we suppose $n\ge 18$. By Theorem \ref{12t}, $H$ and $G$ must be co-degree. Obviously, $H=K_1\vee (sP_1\cup P_{n_{s+1}}\cup P_{n_{s+2}}\cup\cdots \cup P_{n_t})$ or $H=K_1\vee (sP_1\cup P_{n_{s+1}}\cup P_{n_{s+2}}\cup\cdots \cup P_{n_t}\cup C_{q_1}\cup C_{q_2}\cdots \cup C_{q_k})$ ($2\le n_{s+1}\le n_{s+2} \le\cdots \le n_t, k\ge 1$). But now combining this with Corollary \ref{c-multifan}, $H=G$, i.e. $G$ is $DLS$.  \qed

\par\medskip
\noindent{\bf Proof of Corollary \ref{c1.4}.} 
For convenience, let $G:=K_1\vee (P_{l_1}\cup P_{l_1}\cup\cdots \cup P_{l_t}\cup C_{s_1}\cup C_{s_2}\cup\cdots \cup C_{s_k})$ $(t\ge 1, k\ge 1)$  with $n\ge 18$ vertices and $H$ be $L$-cospectral with $G$.  By Corollary \ref{c-multifan} and Lemma \ref{L1-n}, $\lambda_1(G)=n\notin\{\lambda_1(W_3),\lambda_1(W_5)\}$. Moreover, if  each $s_i$ $(i=1,2,\ldots, k)$ is odd, then $\lambda_2(G)<5$ by  Corollary \ref{c-multifan}.   Then by Theorem \ref{12t}, $H$ and $G$ must be co-degree.  Similar with the proof of Corollary \ref{c1.3}, it is easy to see that $H=G$ by Corollary \ref{c-multifan}, i.e. $G$ is $DLS$.
\qed

\section*{Acknowledgements}
The author is deeply grateful to Prof. Muhuo Liu, Prof. Zoran Stani\'c  and Prof. Jianguo Qian for giving some help.

	\end{document}